\documentstyle[11pt]{article}

\def\b#1{{\bf #1}}
\def\i#1{{\it #1}}

\begin{document}

\title{Degenerated Third Order Linear Recurrences}
\author{Mario Catalani\thanks{
Department of Economics, University of Torino, Italy.
mario.catalani@unito.it}}
\date{}
\maketitle

\section{Introduction}
Let us consider a third order linear recurrence
$$U_n=a_1U_{n-1}+a_2U_{n-2}+a_3U_{n-3},\qquad U_0=u_0,\,
U_1=u_1,\,U_2=u_2.$$
Let us assume that the roots of the characteristic polynomial are real,
simple, but two of them are equal in absolute value:
we call this situation
\i{degenerated}. If $\{\lambda_i\}$ are the roots let
us assume that $\lambda_3>0$, $\lambda_1=-\lambda_3$ and
$\lambda_1<\lambda_2<\lambda_3$. Then the coefficients of the
recurrence are $a_1=\lambda_2,\, a_2=\lambda_3^2,\, a_3=-\lambda_2
\lambda_3^2$.

\noindent
The closed form (generalized Binet's formula) is
$$U_n=c_1(-\lambda_3)^n+c_2\lambda_2^n+c_3\lambda_3^n,$$
where $\{c_i\}$ are determined by the initial condition $\{u_i\}$, as well as
by the roots $\{\lambda_i\}$.

\section{Limits}
Using Binet's formula the ratio of two consecutive terms is given by
$${U_n\over U_{n-1}}={c_1(-\lambda_3)^n+c_2\lambda_2^n+c_3\lambda_3^n
\over c_1(-\lambda_3)^{n-1}+c_2\lambda_2^{n-1}+c_3\lambda_3^{n-1}}.$$
After dividing both numerator and denominator by $\lambda_3^{n-1}$ we get
$$
{U_n\over U_{n-1}}={c_1
\lambda_3(-1)^n+c_2\lambda_2\left ({\lambda_2\over\lambda_3}\right )^{n-1}
+c_3\lambda_3\over
c_1(-1)^{n-1}+c_2\left ({\lambda_2\over\lambda_3}\right )^{n-1}
+c_3}.$$
Now the limit for $n\longrightarrow \infty$ depends on the parity of $n$.
If $n$ is even
\begin{equation}
L_1=\lim_{n\longrightarrow\infty}
{U_n\over U_{n-1}}=\lambda_3{c_1+c_3\over-c_1+c_3},
\end{equation}
while if $n$ is odd
\begin{equation}
L_2=\lim_{n\longrightarrow\infty}
{U_n\over U_{n-1}}=\lambda_3{-c_1+c_3\over c_1+c_3}.
\end{equation}
Let us express $c_1$ and $c_3$ explicitly as functions of the initial
conditions. Using Binet's formula we have
$$\left \{\begin{array}{lll}
c_1+c_2+c_3&=&u_0,\\
-c_1\lambda_3+c_2\lambda_2+c_3\lambda_3&=&u_1,\\
c_1\lambda_3^2+c_2\lambda_2^2+c_3\lambda_3^2&=&u_2. \end{array}\right .$$
Let \b{A} be the following matrix
$$\b{A}=\left [\begin{array}{ccc}
1&1&1\\-\lambda_3&\lambda_2&\lambda_3\\
\lambda_3^2&\lambda_2^2&\lambda_3^2 \end{array}\right ].$$
Let $\b{c}$ the column vector with elements $\{c_i\}$
and $\b{u}$ the column vector with elements $\{u_i\}$. Then
$$\b{c}=\b{A}^{-1}\b{u}.$$
We obtain
$$\b{A}^{-1}=\left [\begin{array}{ccc}
{\lambda_2\over 2(\lambda_2+\lambda_3)}&-{1\over 2\lambda_3}&
{1\over 2\lambda_3^2+2\lambda_2\lambda_3}\\
{\lambda_3^2\over\lambda_3^2-\lambda_2^2}&0&{
1\over -\lambda_3^2+\lambda_2^2}\\
{\lambda_2\over -2\lambda_3+2\lambda_2}&{1\over 2\lambda_3}&
{1\over 2\lambda_3^2-2\lambda_2\lambda_3}\end{array}\right ].$$
Then
$$c_1+c_3={-\lambda_2^2u_0+u_2\over \lambda_3^2-\lambda_2^2},$$
and
$$c_3-c_1={\lambda_3^2(-\lambda_2 u_0+u_1) + \lambda_2
(-\lambda_2u_1+u_2)\over \lambda_3^3-\lambda_3\lambda_2^2}.$$
It follows
$$L_1=
{\lambda_3^2(\lambda_2^2 u_0-u_2)\over
\lambda_3^2(\lambda_2 u_0-u_1) + \lambda_2
(\lambda_2u_1-u_2)},$$
$$L_2=
{\lambda_3^2(\lambda_2 u_0-u_1) + \lambda_2
(\lambda_2 u_1-u_2)\over \lambda_2^2u_0-u_2}.$$

\section{Consequences}
If we write
\begin{equation}
\gamma= {c_1+c_3\over-c_1+c_3},
\end{equation}
then
$$
{L_1\over \lambda_3}=\gamma, \qquad {L_2\over \lambda_3}=
{1\over\gamma}.
$$
It follows
\begin{eqnarray}
\lambda_3^2&=&L_1L_2\nonumber\\
&=&
\lim_{n\longrightarrow\infty}
{U_n\over U_{n-1}}
\lim_{n\longrightarrow\infty}
{U_{n-1}\over U_{n-2}}\nonumber\\
&=&
\lim_{n\longrightarrow\infty}
{U_n\over U_{n-1}}{U_{n-1}\over U_{n-2}}\nonumber\\
&=&
\lim_{n\longrightarrow\infty}
{U_n\over U_{n-2}}.
\end{eqnarray}
That is, this limit exists and it does not depend on the initial conditions.

\noindent
Also
\begin{eqnarray}
\gamma^2&=&{L_1\over L_2}\nonumber\\
&=&
{\lim_{n\longrightarrow\infty}
{U_n\over U_{n-1}}\over
\lim_{n\longrightarrow\infty}
{U_{n-1}\over U_{n-2}}}\nonumber\\
&=&
\lim_{n\longrightarrow\infty}
{U_n\over U_{n-1}}{U_{n-2}\over U_{n-1}}\nonumber\\
&=&
\lim_{n\longrightarrow\infty}
{U_nU_{n-2}\over U_{n-1}^2}.
\end{eqnarray}
That is, this limit exists and depends on the initial conditions.

\section{Existence of a Limit}
Under particular initial conditions our ratio converges. The equality
$L_1=L_2$ entails
$$\lambda_3^2(\lambda_2^2u_0-u_2)^2=
[\lambda_3^2(\lambda_2u_0-u_1)+\lambda_2(\lambda_2u_1-u_2)]^2.$$
We can obtain $u_2$ as a function of $u_0$ and $u_1$. This leads to
a quadratic equation in $u_2$ whose solutions are
$$u_2=\left\{\begin{array}{c}\lambda_2\lambda_3u_0+(\lambda_2-
\lambda_3)u_1\\
-\lambda_2\lambda_3u_0+(\lambda_2+
\lambda_3)u_1.\end{array}\right .$$
If we take the first one then
$$L_1=L_2=-\lambda_3.$$
On the other hand
if we take the second one then
$$L_1=L_2=\lambda_3.$$

\end{document}